\documentclass[letterpaper, 10 pt, conference]{ieeeconf}  %

\IEEEoverridecommandlockouts                              %
                                                          
\title{FanCyqlone: A Parallelized Solver for Multi-Horizon Optimal Control}
\author{
  Matisse Teuwen\and
  Pieter Pas\and
  Mathijs Schuurmans\and 
  Panagiotis Patrinos \thanks{This research has received funding from the European Union under the Horizon Europe Grant Agreement AIGGREGATE, no.\;101202457.
  The work was supported by FWO PhD fellowships 11B5027N and 11M9525N and FWO postdoctoral fellowship 1248626N.
	All authors are with the Department of Electrical Engineering (\textsc{esat-stadius}), KU Leuven
	Kasteelpark Arenberg 10, 3001 Leuven, Belgium.
	Email: \texttt{matisse.teuwen@esat.kuleuven.be}}
}

\usepackage{preamble}
\newcommand{\treeQP}{treeQP}
\newcommand{\PIQP}{\textsc{Piqp}}
\newcommand{\FanCyqlone}{\textsc{FanCyqlone}}
\newcommand{\FanCyQpalm}{\textsc{FanCyQPALM}}
\newcommand{\Cyqlone}{\textsc{Cyqlone}}
\newcommand{\QPALM}{\textsc{QPALM}}
\newcommand{\QPALMOCP}{\textsc{QPALM-OCP}}

\begin{document}
\maketitle

\begin{abstract}
  \providecommand{\textscfix}[1]{{\mdseries \textsc{#1}}}%
Multi-horizon optimal control aims to select optimal control sequences 
simultaneously for several possible future scenarios, ensuring consistency at the first stage.
While this enables advanced \acl{MPC} schemes that explicitly handle prediction uncertainty,
solving these problems tends to be costly.
We introduce \textscfix{FanCyqlone}, a linear solver that exploits parallelism 
across both time steps and scenario branches.
The method builds on \textscfix{Cyqlone} to compute the contribution of each branch
and combines these results through a shared root problem.
We integrate \textscfix{FanCyqlone} into \textscfix{QPALM}
to develop {\mdseries \FanCyQpalm{}},
a solver for multi-horizon optimal control problems with general linear stage-wise
constraints.
We compare its performance with state-of-the-art solvers on
spring-mass and quadruped control benchmarks. 

\end{abstract}
\acresetall
\section{Introduction}
\Ac{MPC} is a widely adopted control strategy that optimizes the future performance of a system subject to constraints \cite{MPCbook}.
While standard \ac{MPC} formulations consider a single prediction horizon, many practical applications
call for planning under multiple future scenarios.
This need may arise from uncertainty in the system dynamics \cite{lucia2013multi,gros20172stage},
from multimodal predictions of the behavior of other actors in multi-agent systems \cite{multimodal2025,multimodel2025b}
or from the need for a safe backup trajectory \cite{9838190,geurts2025contingency}.

In this paper, we consider a multi-horizon setting in which scenario trajectories are coupled only at the first stage, 
through a shared root node or through non-anticipativity constraints that enforce identical first-stage control actions across scenarios.
This setting is also referred to as branch \ac{MPC} with a single shared root node \cite{song2026accelerating},
or as multistage robust or scenario \ac{MPC} with a robust horizon of one \cite{lucia2013multi,Klintberg2016dual}.
A prominent special case of multi-horizon \ac{MPC} is contingency \ac{MPC}, which considers two prediction horizons \cite{alsterda2019contingency,alsterda2021contingency}.

Despite the growing interest in multi-horizon and contingency \ac{MPC}, 
the core challenge impeding its widespread adoption is its computational cost, which, 
for a naive implementation may indeed be larger than its deterministic counterpart 
by a factor approximately equal to the number of scenarios.

However, as performance improvements of modern \acp{CPU} are increasingly 
attributable to increasing parallel processing capabilities, new 
opportunities arise to accelerate solution methods for multi-horizon \acp{OCP}.

Motivated by this observation, we propose \FanCyqlone{}, a highly efficient solver for
linear systems arising in multi-horizon linear-quadratic optimal control problems.
It is built on \Cyqlone{} \cite{pas2025cyqlone},
a parallel solver for linear systems
with optimal control structure that exploits the capabilities of
modern \ac{CPU}
architectures through parallelization and batch-wise
vectorization along the time dimension. 
\FanCyqlone{} extends this by additionally parallelizing computations across horizons.
We further integrate \FanCyqlone{} as a linear solver within the
proximal augmented Lagrangian solver \QPALM{} \cite{hermans2019qpalm,hermans2022qpalm,lowenstein2024qpalm},
yielding \FanCyQpalm{}, a solver for constrained multi-horizon linear quadratic optimal control problems.

\subsection{Related work}

In \cite{gros20172stage}, a parallelizable primal-dual interior-point method 
was proposed for two-stage robust \ac{MPC},
which is equivalent to the multi-horizon setting considered here.

The multi-horizon linear optimal control formulation considered in this work is a special case of tree-structured \acp{QP}, which have been studied extensively:
\cite{frison2017high} proposes an interior-point method based on a modified Riccati recursion 
tailored to tree-structured linear \ac{OCP}s,
which was later incorporated into the open-source software package \textsc{Hpipm} \cite{frison2020hpipm}.
However, the method does not parallelize computations across the different scenarios in the tree.

Another line of work decomposes the tree-structured \ac{QP} into smaller, decoupled subproblems that can be solved efficiently in parallel.
This is achieved by dualizing the non-anticipativity constraints that couple the scenarios \cite{Leidereiter2015dual,Klintberg2016dual}, and applying Newton's method to the resulting dual problem.
Similarly, by dualizing both the non-anticipativity and dynamics constraints,
one may obtain smaller stage-wise \acp{QP} at the expense of a larger Newton system \cite{kouzoupis2018dual}.
By only dualizing the dynamics constraints, neighboring scenarios can share decision variables, decreasing the overall number of variables \cite{kouzoupis2019dual}.
Both \cite{kouzoupis2018dual} and \cite{kouzoupis2019dual} are implemented in the open-source library \treeQP{}.

These approaches are particularly effective when the resulting stage-wise \acp{QP} can be solved very efficiently, for example,
in the case of box constraints and diagonal cost Hessians \cite{kouzoupis2018dual,kouzoupis2019dual}.
For tree-structured \acp{QP} with general linear constraints, however, these subproblems can be considerably more expensive to solve. Since the globalization procedure requires solving them repeatedly, their cumulative computational cost can become a limiting factor.

In \cite{schwan2025exploiting}, a proximal interior-point method is proposed for multistage optimization problems with global variables, a problem class that accommodates our formulation, while exploiting the resulting block-tridiagonal arrow structure.
Building on this work, \cite{song2025parallel} develops parallel factorization as well as forward 
and backward solve schemes to further improve computational efficiency.

Recently, there has been growing interest in utilizing the 
massive parallelization capabilities afforded by \acp{GPU}
for solving tree-structured problems.
First-order methods have been used to allow massive parallelization 
of computations in stochastic \cite{sampathirao2024massively} and risk-averse 
\cite{moran2025gpu,bodard2023spock} optimal control problems over scenario trees.
\cite{zhang2025parallel} exploits parallelism both across branches and 
along the time horizon, using a parallel-scan approach to solve the resulting
linear-quadratic subproblems. 
In \cite{sousa2026rake}, a parallel Riccati solver for tree-structured \ac{LQR} problems is developed using rake-compress tree contraction, with \ac{CPU} and \ac{GPU} support.
These \ac{GPU}-based approaches are often particularly effective
for problems involving huge numbers of scenarios. In practical
applications, such as contingency \ac{MPC}, however, the number of
scenarios is often quite modest. In this regime, \ac{GPU} overheads,
including memory-access latency, typically still outweigh the benefits \cite{zhang2025parallel}.
This motivates our focus on \ac{CPU}-based parallelization.

More recently, and concurrently with the present work, \cite{song2026accelerating} proposed a linear solver that exploits the multi-horizon structure by combining parallel Cholesky factorizations with parallel forward and backward substitutions,
thereby achieving two levels of \ac{GPU} parallelism.
However, their approach requires symmetric positive definite system matrices.
As a result, it cannot be directly applied to equality-constrained \ac{KKT} systems,
arising in the inner Newton steps of solvers such as \textsc{Hpipm} and \QPALMOCP{},
as these exhibit a symmetric indefinite saddle-point structure.

\subsection{Contributions}
We highlight the following contributions of this paper.
\begin{enumerate}
\item We develop \FanCyqlone{}, a parallel linear solver on CPU for
    multi-horizon \ac{KKT} systems that extracts branch contributions
    directly from \Cyqlone{}'s parallel-in-time factorization
    and forward solve, applied in parallel across branches. These
    contributions are combined in a single reduced root system,
    after which the branch solutions are recovered through parallel
    backward substitution.
    \item We integrate \FanCyqlone{} into the proximal
    augmented Lagrangian solver \QPALM{}, yielding
    \FanCyQpalm{}, a solver for multi-horizon linear
    optimal control problems with general linear constraints.
    \item We provide an efficient \ac{SPMD} implementation of \FanCyQpalm{}
      in $\Cpp$. Together with Python bindings, it is
      available as open-source software at \href{https://github.com/kul-optec/fancyqlone}{github.com/kul-optec/fancyqlone}.
    \item We demonstrate the computational performance
    of \FanCyQpalm{} on multi-horizon spring-mass
    and quadruped \ac{MPC} benchmarks, showing substantial
    speedups over existing state-of-the-art solvers.
\end{enumerate}
\section*{Notation}
Let $[K]$ denote the set of integers $\{1,2,\ldots,K\}$ and $[K]_0$ denote the set of integers $\{0,1,\ldots,K-1\}$.
The sets of real symmetric positive semidefinite and positive definite matrices of size $n \times n$ are denoted by $\symp{n}$ and $\sympp{n}$, respectively.

\section{Multi-horizon \ac{LQR} Problem}\label{sec:lqr}
We start by considering the following discrete-time multi-horizon \ac{LQR} problem:
\begin{subequations}\label{eq:multi-horz}
    \begin{align}
        &\minimize_{\ivar{x}{\bidx}{\tidx}, \ivar{u}{}{0}, \ivar{u}{\bidx}{\tidx}} && \sum_{\bidx=1}^\numbranch \sum_{\tidx=0}^{N-1} \ivar{\ell}{\bidx}{\tidx}(\ivar{x}{\bidx}{\tidx}, \ivar{u}{\bidx}{\tidx}) + \ivar{\ell}{\bidx}{N}(\ivar{x}{\bidx}{N}) \\
        &\sttlong && \ivar{x}{\bidx}{0} = \ivar{x}{\bidx}{\text{init}},\quad \ivar{u}{\bidx}{0} = \ivar{u}{}{0} \label{eq:initial_state}\\
        & && \ivar{x}{\bidx}{\tidx+1} = \ivar{A}{\bidx}{\tidx} \ivar{x}{\bidx}{\tidx} + \ivar{B}{\bidx}{\tidx} \ivar{u}{\bidx}{\tidx}+ \ivar{f}{\bidx}{\tidx},\\% \quad \forall \tidx \in \{0, \ldots, N-1\}, \forall \bidx \in \{1, \ldots, K\}.
        & && \tidx\in[N]_0,\, j\in[\numbranch]
    \end{align}
\end{subequations}
where the strongly convex stage-wise and terminal costs for each branch are:\vspace{-1.5em}
\[\ivar{\ell}{\bidx}{\tidx}(x,u) = \frac{1}{2} 
\begin{pmatrix}
    u\\x
\end{pmatrix}^\top
\overbrace{
\begin{pmatrix}
\ivar{R}{\bidx}{\tidx}&\ivar{S}{\bidx}{\tidx}\\
\ivar{S}{\bidx}{\tidx\top}&\ivar{Q}{\bidx}{\tidx}
\end{pmatrix}
}^{\displaystyle \ivar{H}{\bidx}{\tidx}}
\begin{pmatrix}
    u\\x
\end{pmatrix}
+
\begin{pmatrix}
    u\\x
\end{pmatrix}^\top
\begin{pmatrix}
    \ivar{r}{\bidx}{\tidx}\\\ivar{q}{\bidx}{\tidx}
\end{pmatrix}
\]
and \(\ivar{\ell}{\bidx}{N}(x) = \frac{1}{2} x^\top \ivar{H}{\bidx}{N} x + x^\top\ivar{q}{\bidx}{N} \) with \(\ivar{H}{\bidx}{N}\coloneq\ivar{Q}{\bidx}{N}\), respectively,
where \(\ivar{Q}{\bidx}{\tidx} \in \sympp{n_x},\ivar{S}{\bidx}{\tidx} \in \R^{n_u \times n_x}\) and \(\ivar{R}{\bidx}{\tidx}\in \sympp{n_u}\) and $\ivar{H}{\bidx}{\tidx}\in\sympp{n_x+n_u}$.
The initial state for each branch is given by \(\ivar{x}{\bidx}{\text{init}}\in \R^{n_x}\).

In \cref{sec:constrained}, we discuss how our method can be used to additionally handle 
linear inequality constraints on the states and inputs.

\subsection{Dynamic programming perspective}
By using backward dynamic programming,
problem \eqref{eq:multi-horz} can be equivalently stated as
\begin{subequations}\label{eq:dp}
    \begin{align}
        &\minimize_{\ivar{x}{\bidx}{0}, \ivar{u}{}{0}, \ivar{x}{\bidx}{1}} && \sum_{\bidx=1}^\numbranch \ivar{\ell}{\bidx}{0}(\ivar{x}{\bidx}{0}, \ivar{u}{}{0}) + \ivar{V}{\bidx}{}(\ivar{x}{\bidx}{1}) \\
        &\sttlong && \ivar{x}{\bidx}{0} = \ivar{x}{\bidx}{\text{init}}, &&& \bidx\in[\numbranch] \\
        & && \ivar{x}{\bidx}{1} = \ivar{A}{\bidx}{0} \ivar{x}{\bidx}{0} + \ivar{B}{\bidx}{0} \ivar{u}{}{0}+ \ivar{f}{\bidx}{0}, &&&\bidx\in[\numbranch] \label{dynamics}%
      \end{align}
\end{subequations}
where
\begin{equation}\label{eq:lqr}
    \begin{aligned}
        \ivar{V}{\bidx}{}(x) = &\min_{\ivar{x}{\bidx}{\tidx}, \ivar{u}{\bidx}{\tidx}} && \sum_{\tidx=1}^{N-1} \ivar{\ell}{\bidx}{\tidx}(\ivar{x}{\bidx}{\tidx}, \ivar{u}{\bidx}{\tidx}) + \ivar{\ell}{\bidx}{N}(\ivar{x}{\bidx}{N}) \\
        & \sttshort && \ivar{x}{\bidx}{1} = x,\\
        & && \ivar{x}{\bidx}{\tidx+1} = \ivar{A}{\bidx}{\tidx} \ivar{x}{\bidx}{\tidx} + \ivar{B}{\bidx}{\tidx} \ivar{u}{\bidx}{\tidx}+ \ivar{f}{\bidx}{\tidx},\\
        & && \tidx\in[N-1] %
      \end{aligned}
\end{equation}
is the cost-to-go function for branch \(\bidx\in[\numbranch]\).
It is well known that the cost-to-go function is a quadratic function of \(x\) \cite[\S 1.3.3]{MPCbook},
i.e.,
\begin{equation}\label{eq:quad-cost-to-go}
\ivar{V}{\bidx}{}(x) = \frac{1}{2} x^\top \ivar{P}{\bidx}{} x + \ivar{p}{\bidx}{\top} x + \ivar{c}{\bidx}{},
\end{equation}
for some $\ivar{P}{\bidx}{} \in \sympp{n_x}$, $\ivar{p}{\bidx}{}\in \R^{n_x}$, and $\ivar{c}{\bidx}{} \in \R$.

We refer to \eqref{eq:dp} as the \textit{root problem}.
In \cref{sec:root,sec:cyqlone}, we 
describe how we can use this decomposition of the problem to efficiently 
obtain a solution to \eqref{eq:multi-horz}.

\subsection{Solving the root problem} \label{sec:root}

After eliminating the initial state for each branch
and using \eqref{eq:quad-cost-to-go}, 
the \ac{KKT} conditions of problem \eqref{eq:dp} are 
\begin{subequations}\label{eq:KKT_dp}
\begin{empheq}[left=\empheqlbrace]{align}
\ivar{P}{\bidx}{} \ivar{x}{\bidx}{1}
- \ivar{\lambda}{\bidx}{0}
&= -\ivar{p}{\bidx}{},
&& \forall j\in[\numbranch]
\label{eq:KKT_dp_a}\\
\ivar{B}{\bidx}{0} \ivar{u}{}{0}
- \ivar{x}{\bidx}{1}
&= -\ivar{\tilde{f}}{\bidx}{0},
&& \forall j\in[\numbranch]
\label{eq:KKT_dp_b}\\
\ivar{R}{}{0} \ivar{u}{}{0}
+ \textstyle\sum_{\bidx=1}^{\numbranch}
\ivar{B}{\bidx}{0\top}\ivar{\lambda}{\bidx}{0}
&= -\tilde{r}^0,
\label{eq:KKT_dp_c}
\end{empheq}
\end{subequations}
with $R^0=\sum_{\bidx=1}^\numbranch \ivar{R}{\bidx}{0}$, $\tilde{r}^0 = \sum_{\bidx=1}^\numbranch \ivar{r}{\bidx}{0}+\ivar{S}{\bidx}{0}\ivar{x}{\bidx}{\text{init}}$, and $\ivar{\tilde{f}}{\bidx}{0} = \ivar{f}{\bidx}{0} + \ivar{A}{\bidx}{0}\ivar{x}{\bidx}{\text{init}}$. The Lagrange multipliers \(\ivar{\lambda}{\bidx}{0}\) are associated with the dynamics constraints \eqref{dynamics} per branch $j$.
Further eliminating $\ivar{x}{\bidx}{1}$ for all $j \in [J]$,
using equations \eqref{eq:KKT_dp_a}, the resulting system can be written in matrix form as:
\begin{equation}\label{eq:KKT_system}
{
\setlength{\arraycolsep}{1pt}
\renewcommand{\arraystretch}{1.15}
\begin{pNiceArray}{ccc|c}
    -P_1^{-\!1} & & & \ivar{B}{1}{0} \\
    & \ddots & & \vdots \\
    & & -P_J^{-\!1} & \ivar{B}{J}{0} \\
    \hline
    \ivar{B}{1}{0\top}
    & \hdots
    & \ivar{B}{J}{0\top}
    & \ivar{R}{}{0}
\end{pNiceArray}
\begin{pmatrix}
    \ivar{\lambda}{1}{0}\\
    \vdots\\
    \ivar{\lambda}{J}{0}\\
    \ivar{u}{}{0}
\end{pmatrix}
=
\begin{pmatrix}
    -\ivar{\tilde{f}}{1}{0}-P_1^{-\!1}\ivar{p}{1}{}\\
    \vdots\\
    -\ivar{\tilde{f}}{J}{0}-P_J^{-\!1}\ivar{p}{J}{}\\
    -\tilde{r}^0
\end{pmatrix}}
\end{equation}
Following the block Cholesky factorization procedure as outlined in \cite[\S 3.1]{pas2025cyqlone}, we obtain the following factorization procedure of the system matrix:
\begin{Fsteps}   
     \item Factorize each $\ivar{P}{\bidx}{-\!1}=\ivar{L}{\bidx}{P} \ivar{L}{\bidx}{P\top}$ (parallel);\label{F1}
     \item Calculate each \(\ivar{L}{\bidx}{B} = \ivar{B}{\bidx}{0\top}\ivar{L}{\bidx}{P-\!\top}\) (parallel);
     \item Calculate each \(\ivar{S}{\bidx}{B} = \ivar{L}{\bidx}{B}\ivar{L}{\bidx}{B\top}\) (parallel);
     \item Calculate the Schur complement \(\mathscr{R}=R^0 + \sum_{\bidx=1}^\numbranch \ivar{S}{\bidx}{B}\);
     \item Calculate the Cholesky factorization of \(\mathscr{R}=L^{\mathscr R}L^{\mathscr R\top}\).
\end{Fsteps}
In terms of these matrices, the LDL\(^\top\)-factorization of the system matrix is then given by
\begin{equation}\label{eq:rootfactor}{
\setlength{\arraycolsep}{1pt}
\renewcommand{\arraystretch}{1.15}
L=\begin{pNiceArray}{ccc|c}[columns-width=1em ]
    \ivar{L}{1}{P}              &             &             &   \\
                              & \ddots &             & \\
                             &             & \ivar{L}{\numbranch}{P} &  \\
    \hline
         -\ivar{L}{1}{B}
        &  \hdots
        &  -\ivar{L}{\numbranch}{B}
        & L^{\mathscr R} \\ 
\end{pNiceArray},}
\quad{
\setlength{\arraycolsep}{1pt}
\renewcommand{\arraystretch}{1.15}
 D=\begin{pNiceArray}{ccc|c}[columns-width=1em ]
    -\mathrm{I}&             &                          &  \\
      &    \ddots   &                    & \\
      &             &              -\mathrm{I}          &  \\    \hline
      &             &                         & \mathrm{I}\\
\end{pNiceArray}.
}
\end{equation}
The triangular solve procedure for the system \eqref{eq:KKT_system} can be performed in a similar manner, where the forward and backward substitutions are performed in parallel across branches.
That is,
\begin{Ssteps}
    \item Solve each \(\ivar{L}{\bidx}{P}\ivar{\hat{p}}{\bidx}{0}=\ivar{\tilde{f}}{\bidx}{0}+P_\bidx^{-\!1}\ivar{p}{\bidx}{}\)(parallel);\label{S1}
    \item Calculate each $\ivar{s}{\bidx}{} = \ivar{L}{\bidx}{B}\ivar{\hat{p}}{\bidx}{0}$ (parallel);
    \item Solve \(L^{\mathscr R}L^{\mathscr R\top}\ivar{u}{}{0}=-\tilde{r}^0-\sum_{\bidx=1}^{\numbranch}\ivar{s}{\bidx}{}\);
    \item Calculate each \(\ivar{y}{\bidx}{0} = \ivar{\hat{p}}{\bidx}{0} +\ivar{L}{\bidx}{B\top} \ivar{u}{}{0}\) (parallel);
    \item Solve each \(\ivar{L}{\bidx}{P\top}\ivar{\lambda}{\bidx}{0}=\ivar{y}{\bidx}{0}\) (parallel).
\end{Ssteps}
\subsection{Solving the branch problems with \texorpdfstring{\Cyqlone}{Cyqlone}} \label{sec:cyqlone}
The root factorization and solve require solving problem \eqref{eq:lqr} for each branch, 
to obtain $P_j$ and $p_j$.
This could be done using the sequential Riccati recursion,
which is a well-known method for solving \ac{LQR} problems.
However, to enable further parallelization across the time axis,
we instead propose to use \Cyqlone{} \cite{pas2025cyqlone} to solve \eqref{eq:lqr} for all $\bidx \in [\numbranch]$ (in parallel).
In fact, by augmenting the problem with a dummy initial state and input,
$\ivar{L}{\bidx}{P}$ and $\ivar{\hat{p}}{\bidx}{0}$ can be obtained directly 
from the output of \Cyqlone{}, as we now describe.

\subsubsection{Problem augmentation}
\newcommand\Perm[1]{\mathrm P_{\!#1}}
Consider the \ac{LQR} problem \eqref{eq:lqr} for each branch $\bidx \in [\numbranch]$.
We prepend a dummy state $\ivar{\tilde{x}}{\bidx}{0}$ and a dummy input $\ivar{\tilde{u}}{\bidx}{0}$, with dynamics
\(
\ivar{x}{\bidx}{1}=\ivar{\tilde{x}}{\bidx}{0},
\)
and the corresponding stage cost is chosen as
\(
\ivar{\tilde{\ell}}{\bidx}{0}(x,u)
=\frac{1}{2}u^\top u\).
We refer to the resulting \ac{LQR} problem:
\begin{equation}\label{eq:auglqr}
    \begin{aligned}
        &\minimize_{\ivar{\tilde{x}}{\bidx}{0},\ivar{\tilde{u}}{\bidx}{0},\ivar{x}{\bidx}{\tidx}, \ivar{u}{\bidx}{\tidx}} && \ivar{\tilde{\ell}}{\bidx}{0}(\ivar{\tilde{x}}{\bidx}{0},\ivar{\tilde{u}}{\bidx}{0}) + \sum_{\tidx=1}^{N-1} \ivar{\ell}{\bidx}{\tidx}(\ivar{x}{\bidx}{\tidx}, \ivar{u}{\bidx}{\tidx}) + \ivar{\ell}{\bidx}{N}(\ivar{x}{\bidx}{N})\\
        & \sttshort && \ivar{x}{\bidx}{1}= \ivar{\tilde{x}}{\bidx}{0}=  x,\\
        & && \ivar{x}{\bidx}{\tidx+1} = \ivar{A}{\bidx}{\tidx} \ivar{x}{\bidx}{\tidx} + \ivar{B}{\bidx}{\tidx} \ivar{u}{\bidx}{\tidx}+ \ivar{f}{\bidx}{\tidx},\\
        & && \tidx\in[N-1], %
      \end{aligned}
\end{equation}
as the \textit{augmented} LQR problem.
\subsubsection{\Cyqlone{} algorithm}

\begin{wrapfigure}{r}{0.33\linewidth}
    \centering
    \begin{tikzpicture}[scale= 0.30]
\draw[fill={lightgray},
draw={white}](0.000000000, 0.000000000) rectangle (1.000000000, -1.000000000);
\draw[fill={gray},
draw={white}](8.000000000, 0.000000000) rectangle (8.700000000, -1.000000000);
\draw[fill={lightgray},
draw={white}](1.000000000, -1.000000000) rectangle (2.000000000, -2.000000000);
\draw[fill={gray},
draw={white}](8.000000000, -1.000000000) rectangle (8.700000000, -2.000000000);
\draw[fill={lightgray},
draw={white}](2.000000000, -2.000000000) rectangle (3.000000000, -3.000000000);
\draw[fill={gray},
draw={white}](8.000000000, -2.000000000) rectangle (8.700000000, -3.000000000);
\draw[fill={lightgray},
draw={white}](3.000000000, -3.000000000) rectangle (4.000000000, -4.000000000);
\draw[fill={gray},
draw={white}](8.000000000, -3.000000000) rectangle (8.700000000, -4.000000000);
\draw[fill={lightgray},
draw={white}](4.000000000, -4.000000000) rectangle (5.000000000, -5.000000000);
\draw[fill={gray},
draw={white}](8.000000000, -4.000000000) rectangle (8.700000000, -5.000000000);
\draw[fill={lightgray},
draw={white}](5.000000000, -5.000000000) rectangle (6.000000000, -6.000000000);
\draw[fill={gray},
draw={white}](8.000000000, -5.000000000) rectangle (8.700000000, -6.000000000);
\draw[fill={lightgray},
draw={white}](6.000000000, -6.000000000) rectangle (7.000000000, -7.000000000);
\draw[fill={gray},
draw={white}](8.000000000, -6.000000000) rectangle (8.700000000, -7.000000000);
\draw[fill={lightgray},
draw={white}](7.000000000, -7.000000000) rectangle (8.000000000, -8.000000000);
\draw[fill={gray},
draw={white}](8.000000000, -7.000000000) rectangle (8.700000000, -8.000000000);
\draw[fill={gray},
draw={white}](0.000000000, -8.000000000) rectangle (1.000000000, -8.700000000);
\draw[fill={gray},
draw={white}](1.000000000, -8.000000000) rectangle (2.000000000, -8.700000000);
\draw[fill={gray},
draw={white}](2.000000000, -8.000000000) rectangle (3.000000000, -8.700000000);
\draw[fill={gray},
draw={white}](3.000000000, -8.000000000) rectangle (4.000000000, -8.700000000);
\draw[fill={gray},
draw={white}](4.000000000, -8.000000000) rectangle (5.000000000, -8.700000000);
\draw[fill={gray},
draw={white}](5.000000000, -8.000000000) rectangle (6.000000000, -8.700000000);
\draw[fill={gray},
draw={white}](6.000000000, -8.000000000) rectangle (7.000000000, -8.700000000);
\draw[fill={gray},
draw={white}](7.000000000, -8.000000000) rectangle (8.000000000, -8.700000000);
\draw[fill={black},
draw={white}](8.000000000, -8.000000000) rectangle (8.700000000, -8.700000000);
\draw[decorate,
decoration={brace, amplitude=5pt}](1.000000000, 0.000000000) -- (8.000000000, -7.000000000);
\node[anchor={south},
rotate={-45.0},
yshift={5pt}] at (4.500000000, -3.500000000){$\np \nv$ blocks};

\end{tikzpicture}
    \vspace{5pt}
    \caption{Block structure of the permuted KKT matrix
        $\mathscr{K}_j$.}
    \label{fig:cyqlone-kkt}
\end{wrapfigure}
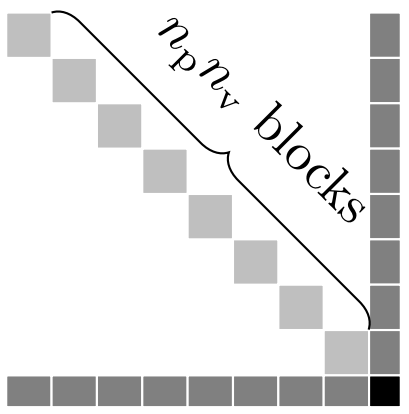

\Cyqlone{} is a linear solver specialized for \ac{KKT} systems with optimal control 
structure.
We use it to solve a problem of the form \eqref{eq:auglqr} for each branch $j$.
After eliminating the initial state $\ivar{\tilde{x}}{\bidx}{0}$,
\Cyqlone{} permutes the corresponding \ac{KKT} system
as described in \cite[\S 4.1]{pas2025cyqlone}, resulting in a 
matrix $\mathscr{K}_\bidx$ whose structure is shown in \Cref{fig:cyqlone-kkt}.
The banded matrices along the light gray block-diagonal can be factorized independently, using a modified Riccati recursion.
Each of the \(\np\) processors factorizes \(\nv\) blocks
at once  using batchwise \ac{SIMD} vectorization
\cite[\S 7.3]{pas2025cyqlone}.
Eliminating the light gray blocks introduces fill-in in the bottom-right block, shown in black.
The resulting Schur complement is block-tridiagonal and can thus be factorized in parallel using cyclic reduction. We refer to the full factorization procedure as \textsc{CyqloneFact}, and label the resulting Cholesky factor $L^{\mathscr{K}}_\bidx$.
The factors are then used to solve the \ac{KKT} system
of problem \eqref{eq:auglqr} via forward and backward substitution, 
each again parallelized across the independent blocks. We refer to these procedures by \textsc{CyqloneFSolve} and \textsc{CyqloneBSolve}, respectively.
Given the right-hand side vectors $\tilde f_j$ (dynamics constraints), $\tilde r_j$ (control gradients) and $q_j$ (state gradients),
the forward solve computes the corresponding intermediate vectors $\hat{\lambda}_\bidx$, $\hat{u}_\bidx$ and $\hat{x}_\bidx$,
which are then used in the backward solve to compute the solution of the branch variables $\lambda_\bidx$ (Lagrange multipliers), $u_\bidx$ (controls) and $x_\bidx$ (states).

\subsubsection{Connection to the root problem}

We can now describe how the results of \Cyqlone{} can be used to obtain $L^P_j$ and $\ivar{\hat{p}}{\bidx}{0}$
described in \cref{sec:root}.
\begin{proposition}\label{prop:factor}
Consider the augmented \ac{LQR} problem \eqref{eq:auglqr}.
Let $\mathscr{K}_\bidx$ denote its \ac{KKT} matrix under \Cyqlone{}'s permutation \cite[\S 4.1]{pas2025cyqlone},
with $L_{0,j}$ the bottom right $n_x\times n_x$-block of the block Cholesky factor $L^{\mathscr{K}}_\bidx$ of $\mathscr{K}_\bidx$ computed by \Cyqlone{}.
Let $\ivar{\hat{\lambda}}{j}{0}$ denote the corresponding value after the forward solve with $L^{\mathscr{K}}_\bidx$.
Then \[L_{0,j}L_{0,j}^\top=P^{-\!1}_j,\qquad L_{0,j}\ivar{\hat{\lambda}}{j}{0}=x+P^{-\!1}_{\bidx}p_{\bidx}.\]
Furthermore, \(L_{0,j}=L^P_j\) and if $x=\ivar{\tilde{f}}{\bidx}{0}$ then \(\ivar{\hat{\lambda}}{j}{0}=\ivar{\hat{p}}{\bidx}{0}\).
\end{proposition}
\begin{proof}
See Appendix.
\end{proof}

\subsection{\FanCyqlone{}}
\FanCyqlone{}'s factorization \cref{alg:fact-fancyqlone} consists of two stages.
First, on line \ref{line:cyfact} \Cyqlone{}'s factorization algorithm is applied to each augmented branch problem \eqref{eq:auglqr}.
In line \ref{line:extract}, the $L_{0,j}$ factors are extracted as the bottom-right $n_x\times n_x$ blocks of the factors $L^{\mathscr{K}}_j$, denoted as $\operatorname{lastdiagblock}$. The root problem is then factorized in lines \ref{line:rootfact}-\ref{line:rootfactend}.

\FanCyqlone{}'s triangular solve procedure is summarized in \cref{alg:solve-fancyqlone}.
Here, line \ref{line:cyfsolve} applies \Cyqlone{}'s forward substitution algorithm to each augmented branch problem.
Lines \ref{line:rootsolve}-\ref{line:rootsolveend} then solve the root problem, and line \ref{line:cybsolve} applies \Cyqlone{}'s backward substitution algorithm to each augmented branch problem to obtain the solution of the remaining branch variables.
\begin{remark}
    \Cyqlone{} eliminates the initial state $\ivar{\tilde{x}}{\bidx}{0}=x$, so $\ivar{A}{\bidx}{0}$, $\ivar{Q}{\bidx}{0}$ and $\ivar{S}{\bidx}{0}$ are not accessed during the factorization,
    which is why we do not introduce $\widetilde{A}_\bidx$, $\widetilde{Q}_{\bidx}$, and $\widetilde{S}_{\bidx}$.
\end{remark}

\begin{algorithm2e}[t]
    \caption{\FanCyqlone{} factorization}
    \label{alg:fact-fancyqlone}
    \DontPrintSemicolon

    \KwIn{\ac{OCP} data
        $(A_\bidx,B_\bidx,Q_\bidx,R_\bidx,S_\bidx)$,
        for each branch $\bidx\in[J]$}
    \KwIn{$N$, $\npj$, $\nv$: horizon length,
        number of processors per branch and vector length}
    \KwOut{$L_\bidx^{\mathscr K}$, $L_\bidx^B$ for each
        branch, and $L^{\mathscr R}$ 
      }
    $R^0=\sum_{\bidx=1}^{\numbranch}\ivar{R}{\bidx}{0}$\nllabel{line:rootsum}\;
        \tcp{Independent branch factorizations}
    \For({\tcp*[f]{in parallel}}){$\bidx=1,\ldots,\numbranch$}{

    \tcp{{Augment the \ac{LQR} problem }}
    
        $\widetilde B_\bidx \gets
        (0,B_\bidx^1,\ldots,B_\bidx^{N-1})$ \nllabel{line:augmentB}\;

        $\widetilde R_\bidx \gets
        (\mathrm I,R_\bidx^1,\ldots,R_\bidx^{N-1})$ \nllabel{line:augmentR}\; 
        \tcp{Factorize the augmented branch}
        $\begin{aligned}
        L_\bidx^{\mathscr{K}}
        \gets {}&
        \textsc{CyqloneFact}(
             A_\bidx,
            \widetilde B_\bidx,
             Q_\bidx,\\[-0.5 em]
        &\qquad\widetilde R_\bidx,
          S_\bidx,
         N,\npj, \nv)
        \end{aligned}$ \nllabel{line:cyfact}\;
        \tcp{Extract the final Cholesky block}
        $\ivar{L}{0,\bidx}{} \gets
        \operatorname{lastdiagblock}
        (L_\bidx^{\mathscr K})$\nllabel{line:extract}\;

        \tcp{Compute the branch contribution}
        $L_\bidx^B \gets
        B_\bidx^{0\top} \ivar{L}{0, \bidx}{-\!\top}$ \nllabel{line:rootfact}\;

        $S_\bidx^B \gets
        L_\bidx^B L_\bidx^{B\top}$\nllabel{line:endfor}\;
    }
    $\mathscr R \gets
    R^0+\textstyle\sum_{\bidx=1}^{\numbranch}
    S_\bidx^B$\nllabel{line:sum}\;
    \tcp{Factorize the reduced root Hessian}
    $L^{\mathscr R} \gets
    \operatorname{chol}(\mathscr R)$\nllabel{line:rootfactend}\;
\end{algorithm2e}

\begin{algorithm2e}[t]
    \caption{\FanCyqlone{} Triangular Solve}
    \label{alg:solve-fancyqlone}
    \DontPrintSemicolon

    \KwIn{\ac{OCP} data
        $(r_\bidx,q_\bidx,f_\bidx,\ivar{A}{\bidx}{0},\ivar{S}{\bidx}{0},\ivar{x}{\bidx}{\text{init}})$,
        for each branch $\bidx\in[J]$}
    \KwIn{$N$, $\npj$, $\nv$: horizon length,
        number of processors per branch, vector length}
    \KwIn{$L_\bidx^{\mathscr K}$, $L_\bidx^B$ for each
        branch, and $L^{\mathscr R}$}
    \KwOut{$\ivar{\lambda}{\bidx}{\tidx}$, $\ivar{u}{\bidx}{\tidx}$, $\ivar{x}{\bidx}{\tidx}$, $\ivar{u}{}{0}$}
        \tcp{Independent branch forward solve}
    \For({\tcp*[f]{in parallel}}){$\bidx=1,\ldots,\numbranch$}{
        $\ivar{\tilde{f}}{\bidx}{0} = \ivar{f}{\bidx}{0} + \ivar{A}{\bidx}{0}\ivar{x}{\bidx}{\text{init}}$\;

        \tcp{{Augment the \ac{LQR} problem}}
    
        $ \ivar{\tilde r}{\bidx}{} \gets
        (0,\ivar{r}{\bidx}{1},\ldots,\ivar{r}{\bidx}{N-1})$\;
        $ \ivar{\tilde f}{\bidx}{} \gets
        (\ivar{\tilde{f}}{\bidx}{0},\ivar{f}{\bidx}{1},\ldots,\ivar{f}{\bidx}{N-1})$\;
        \tcp{Forward solve the augmented branch}
        $\begin{aligned}
        (\ivar{\hat{\lambda}}{\bidx}{\tidx},\ivar{\hat{u}}{\bidx}{\tidx},\ivar{\hat{x}}{\bidx}{\tidx})
        \gets {}&
        \textsc{CyqloneFSolve}(
             \ivar{L}{\bidx}{\mathscr{K}}, \ivar{\tilde{f}}{\bidx}{}, \ivar{\tilde{r}}{\bidx}{},\\[-0.5 em]
        &\qquad\ivar{q}{\bidx}{}, N,\npj, \nv)
        \end{aligned}$\nllabel{line:cyfsolve}\;        

        \tcp{Compute the branch contribution}
        $\ivar{s}{\bidx}{0} \gets
        \ivar{L}{\bidx}{B}\ivar{\hat{\lambda}}{\bidx}{0}$\nllabel{line:rootsolve}\;
    }
    $\ivar{\tilde{r}}{}{0}=\sum_{\bidx=1}^{\numbranch}\ivar{r}{\bidx}{0}+\ivar{S}{\bidx}{0}\ivar{x}{\bidx}{\text{init}}$\;
    $\ivar{\tilde{r}}{}{0}\gets\ivar{\tilde{r}}{}{0}+\sum_{\bidx=1}^{\numbranch}\ivar{s}{\bidx}{0}$\\
    \tcp{Solve the reduced root system}
    $\ivar{u}{}{0}\gets
    -L^{\mathscr R-\top}L^{\mathscr R-1}\ivar{\tilde{r}}{}{0}$\;
    \For({\tcp*[f]{in parallel}}){$\bidx=1,\ldots,\numbranch$}{
        \tcp{Distribute the root to each branch}
        $\ivar{\hat{\lambda}}{\bidx}{0}\gets
        \ivar{\hat{\lambda}}{\bidx}{0}+\ivar{L}{\bidx}{B\top}\ivar{u}{}{0}$ \nllabel{line:rootsolveend}\;
        \tcp{Backward solve the augmented branch}
        $\begin{aligned}
        (\ivar{\lambda}{\bidx}{\tidx},\ivar{u}{\bidx}{\tidx},\ivar{x}{\bidx}{\tidx})
        \gets {}&
        \textsc{CyqloneBSolve}(
             \ivar{L}{\bidx}{\mathscr{K}}, \ivar{\hat{\lambda}}{\bidx}{\tidx},\ivar{\hat{u}}{\bidx}{\tidx}, \\[-0.5 em]
        &\qquad\ivar{\hat{x}}{\bidx}{\tidx}, N,\npj, \nv)
        \end{aligned}$\nllabel{line:cybsolve}\;        
    }
\end{algorithm2e}

\section{Inequality-Constrained Multi-Horizon \ac{OCP}}\label{sec:constrained}
We consider the \textit{constrained} discrete-time multi-horizon linear-quadratic \ac{OCP}:
\begin{subequations}\label{eq:multi-horz_ineq}
    \begin{align}
        &\minimize_{\ivar{x}{\bidx}{\tidx}, \ivar{u}{}{0}, \ivar{u}{\bidx}{\tidx}} && \sum_{\bidx=1}^\numbranch \sum_{\tidx=0}^{N-1} \ivar{\ell}{\bidx}{\tidx}(\ivar{x}{\bidx}{\tidx}, \ivar{u}{\bidx}{\tidx}) + \ivar{\ell}{\bidx}{N}(\ivar{x}{\bidx}{N}) \\
        &\sttlong && \ivar{x}{\bidx}{0} = \ivar{x}{\bidx}{\text{init}},\quad \ivar{u}{\bidx}{0} = \ivar{u}{}{0} \label{eq:initial_state}\\
        & && \ivar{x}{\bidx}{\tidx+1} = \ivar{A}{\bidx}{\tidx} \ivar{x}{\bidx}{\tidx} + \ivar{B}{\bidx}{\tidx} \ivar{u}{\bidx}{\tidx}+ \ivar{f}{\bidx}{\tidx},\\% \quad \forall \tidx \in \{0, \ldots, N-1\}, \forall \bidx \in \{1, \ldots, K\}.
        & && \ivar{\underline{b}}{\bidx}{\tidx} \leq \ivar{C}{\bidx}{\tidx}\ivar{x}{\bidx}{\tidx} + \ivar{D}{\bidx}{\tidx} \ivar{u}{\bidx}{\tidx} \leq \ivar{\Bar{b}}{\bidx}{\tidx},\label{eq:stage_ineq}\\% \quad \forall \tidx \in \{0, \ldots, N-1\}, \forall \bidx \in \{1, \ldots, K\}.
        & && \ivar{\underline{b}}{\bidx}{N} \leq \ivar{C}{\bidx}{N}\ivar{x}{\bidx}{N} \leq \ivar{\Bar{b}}{\bidx}{N},\label{eq:terminal_ineq}\\% \quad \forall \tidx = N, \forall \bidx \in \{1, \ldots, K\}.    
        & && \tidx\in[N]_0, j\in[\numbranch],
    \end{align}
\end{subequations}
where we now assume only convexity of the costs, i.e., \(\ivar{H}{\bidx}{\tidx}\in\symp{n_x+n_u}\) and \(\ivar{H}{\bidx}{N}\in\symp{n_x}\) 
are positive semidefinite.
\subsection{\FanCyQpalm{}}\label{sec:fan_cyqpalm}
We handle the stage-wise and terminal inequality constraints of the problem using the proximal \ac{ALM} \QPALM{} \cite{hermans2019qpalm,hermans2022qpalm}.
Similar to \QPALMOCP{} \cite{lowenstein2024qpalm} and \textsc{CyQPALM} \cite{pas2025cyqlone}, only the inequality constraints \eqref{eq:stage_ineq} and \eqref{eq:terminal_ineq} are relaxed using \ac{ALM}, while the equality constraints are kept in the inner problem.
Thus, at each outer iteration we minimize the piecewise-quadratic regularized augmented Lagrangian function subject to the equality constraints and update the Lagrange multipliers and penalties.
The inner problems are iteratively solved using semismooth Newton's method. The resulting linear systems have the same structure as the \ac{KKT} system of the multi-horizon \ac{LQR} problem as outlined in \cref{sec:lqr}. They can therefore be solved using \textsc{FanCyqlone}.
The generalized Hessians for these semismooth Newton steps \cite[Eq. 31]{pas2025cyqlone} are given by:
\[\ivar{\widehat{R}}{}{0}=\ivar{R}{}{0}+\textstyle\sum_{\bidx=1}^{\numbranch}
\ivar{D}{\bidx}{0\top} \ivar{\Sigma}{\bidx}{\mathcal{J},0}
\ivar{D}{\bidx}{0}+ \Gamma_u^{-1}\]
\[\ivar{\widehat{H}}{\bidx}{\tidx}=\ivar{H}{\bidx}{\tidx}+\begin{pmatrix}
\ivar{D}{\bidx}{\tidx} & \ivar{C}{\bidx}{\tidx}
\end{pmatrix}^\top \ivar{\Sigma}{\bidx}{\mathcal{J},\tidx}\begin{pmatrix}
\ivar{D}{\bidx}{\tidx} & \ivar{C}{\bidx}{\tidx}
\end{pmatrix}+ \Gamma^{-1}\]
\[\ivar{\widehat{H}}{\bidx}{N}=\ivar{H}{\bidx}{N}+
\ivar{C}{\bidx}{N\top} \ivar{\Sigma}{\bidx}{\mathcal{J},N}
\ivar{C}{\bidx}{N}+ \Gamma_x^{-1}\] with $j\in[J]$ and $\tidx\in[N-1]$.
Here \(\ivar{\Sigma}{\bidx}{\mathcal{J},\tidx}\) is a diagonal penalty matrix with non-zero elements corresponding to the active set $\mathcal{J}$ of the inequality constraints
at the current iterate and $\Gamma= \left( \begin{smallmatrix}
  \Gamma_u & \\
   & \Gamma_x
\end{smallmatrix} \right) $, where $\Gamma_u\in\sympp{n_u}$ and $\Gamma_x\in\sympp{n_x}$ are positive-definite primal regularization matrices ensuring $\ivar{\widehat{H}}{\bidx}{\tidx}\in \sympp{n_x+n_u}$ and $\ivar{\widehat{H}}{\bidx}{N}\in \sympp{n_x}$.

\section{Numerical Experiments}\label{sec:experiments}
The performance of the proposed algorithms is evaluated on two benchmarks. The first is the widely used linear spring-mass benchmark \cite{BoydandWang2010},
which we adapt to the multi-horizon setting by introducing branch-dependent spring constants. 
The second is based on a quadruped centroidal \ac{OCP}
formulation \cite{lezia2022}, where branch-specific friction coefficients define different scenarios.
All benchmarks are performed on an Intel Core Ultra 7 265 \ac{CPU}, with execution restricted to its eight performance cores.

We compare the performance of \FanCyQpalm{} against three state-of-the-art solvers: 
\begin{enumerate}
  \item \textsc{Hpipm} specialized for tree-structured \acp{QP} \cite{frison2017high,frison2020hpipm} using the speed profile%
\footnote{Additionally, we set $\texttt{t0\_init}=0$ for the Quadruped benchmark.
We observed that initializing the slacks and their multipliers in this way reduces numerical issues.};
\item \treeQP{} with the TDunes backend \cite{kouzoupis2019dual}; and 
\item \PIQP{} with its multistage parallel backend \cite{song2025parallel}.
\end{enumerate}
To ensure a fair comparison, for each solver we explicitly eliminate the initial state as in \FanCyQpalm{}. 
Furthermore, except for \FanCyQpalm{} we use bound constraints where possible,
rather than formulating them as general linear constraints, as this decreases solve times considerably.
For \FanCyQpalm{}, $\nv$ denotes the vector length and $\np$ the total number of processors. 
The processors are distributed across the branches, giving $\npj = \frac{\np}{\numbranch}$ processors per branch, assuming $\np\geq\numbranch$.

In order to assess the benefits of parallelizing over both the branches and time,
we consider a \texttt{sequential} variant of \FanCyQpalm{},
in addition to the default \texttt{parallel} variant.
In the \texttt{sequential} variant, the different branches $j \in [J]$ are handled sequentially,
leaving only parallelization over time.
In this case, each branch therefore uses all $\np$ available processors, such that $\npj=\np$.

Since \textsc{Hpipm} only provides a sequential implementation, it is run on a single processor, while \treeQP{} and \PIQP{} are run on $\np$ processors.
We set the OpenMP environment variables \texttt{OMP\_DYNAMIC}=FALSE and \texttt{OMP\_WAIT\_POLICY}=ACTIVE.
We set the absolute convergence tolerances of all solvers to $10^{-8}$. However, to account for differences in their stopping criteria, we disabled \textsc{Piqp}'s duality-gap check and relaxed \textsc{Hpipm}'s complementarity tolerance to $10^{-6}$, as these criteria are not explicitly enforced in the same way by all solvers.
\subsection{Spring-Mass System} \label{sec:spring-mass}

The spring-mass system follows the benchmark problem used in \cite{pas2025cyqlone}.
The system consists of a chain of $M$ masses connected by springs, where each mass can move in one dimension.
For branch $\bidx$, the state vector $\ivar{x}{\bidx}{} \in \mathbb{R}^{2M}$ contains the positions and velocities of the masses. Actuators are placed between pairs of masses, resulting in a number of inputs $n_u=\frac{M}{2}$ for even $M\geq6$. The problem has only box constraints on the inputs and positions.
All benchmark parameters are identical to those used in \cite{pas2025cyqlone},
except for the spring constant $\ivar{k}{\bidx}{}$,
which varies across branches and therefore results in branch-specific system dynamics.
The spring constants are equally spaced between $0.8$ and $1.2$.

\begin{figure}[h!]
    \centering
    \begin{subfigure}{\columnwidth}
        \centering
        \begin{tikzpicture}

\definecolor{scalingcolor0}{HTML}{1f77b4}
\definecolor{scalingcolor1}{HTML}{ff7f0e}
\definecolor{scalingcolor2}{HTML}{2ca02c}
\definecolor{scalingcolor3}{HTML}{d62728}

\begin{axis}[
    xlabel={Horizon length $N$},
    ylabel={Time per iteration [ms]},
    width={\textwidth},
    height={0.6\textwidth},
    ymin={0},
    ymax={3.098687904375},
    grid={major},
    enlarge x limits={0.05},
    enlarge y limits={false},
    scaled ticks={false},
    tick align={outside},
    tick pos={left},
    major grid style={draw=black,draw opacity=0.15},
    legend pos={north west},
    legend cell align={left},
    xtick={32,64,96,128,160,192,224,256},
    legend style={
        font=\fontsize{5.5}{6.5}\selectfont,
        /pgfplots/legend columns=2,
        /pgfplots/transpose legend=false,
        draw=black!20,
        fill=white,
        fill opacity=0.6,
        text opacity=1,
        column sep=0.5em,
        row sep=-1pt,
        inner xsep=2pt,
        inner ysep=1pt,
    },
]

\addplot[
    color={scalingcolor1},
    mark={*},
    mark options={solid},
    line width={0.75pt},
    mark size={0.8pt},
    solid
]
table[row sep={\\}]{%
32.000000000 0.295142739\\
64.000000000 0.559701803\\
96.000000000 0.887656076\\
128.000000000 1.157566623\\
160.000000000 1.475120889\\
192.000000000 1.841557547\\
224.000000000 2.263807941\\
256.000000000 2.733504333\\
};
\addlegendentry{$\np\!=\!4,\ \npj\!=\!1,\ \nv\!=\!1$}

\addplot[
    color={scalingcolor1},
    mark={*},
    mark options={solid},
    line width={0.75pt},
    mark size={0.8pt},
    dashed
]
table[row sep={\\}]{%
32.000000000 0.356112160\\
64.000000000 0.630633805\\
96.000000000 0.972864231\\
128.000000000 1.252317871\\
160.000000000 1.596723407\\
192.000000000 1.958396257\\
224.000000000 2.424601429\\
256.000000000 2.951131338\\
};
\addlegendentry{$\np\!=\!4,\ \npj\!=\!4,\ \nv\!=\!1$}
\addplot[
    color={scalingcolor3},
    mark={*},
    mark options={solid},
    line width={0.75pt},
    mark size={0.8pt},
    solid
]
table[row sep={\\}]{%
32.000000000 0.169867552\\
64.000000000 0.309603934\\
96.000000000 0.478263596\\
128.000000000 0.615751186\\
160.000000000 0.764419621\\
192.000000000 0.917261890\\
224.000000000 1.101486789\\
256.000000000 1.286994176\\
};
\addlegendentry{$\np\!=\!8,\ \npj\!=\!2,\ \nv\!=\!1$}

\addplot[
    color={scalingcolor3},
    mark={*},
    mark options={solid},
    line width={0.75pt},
    mark size={0.8pt},
    dashed
]
table[row sep={\\}]{%
32.000000000 0.245955970\\
64.000000000 0.382644186\\
96.000000000 0.557624487\\
128.000000000 0.698724606\\
160.000000000 0.849199032\\
192.000000000 1.004739189\\
224.000000000 1.203423589\\
256.000000000 1.411803219\\
};
\addlegendentry{$\np\!=\!8,\ \npj\!=\!8,\ \nv\!=\!1$}
\addplot[
    color={scalingcolor0},
    mark={*},
    mark options={solid},
    line width={0.75pt},
    mark size={0.8pt},
    solid
]
table[row sep={\\}]{%
32.000000000 0.111871505\\
64.000000000 0.198500409\\
96.000000000 0.307581298\\
128.000000000 0.417172915\\
160.000000000 0.544732088\\
192.000000000 0.673712707\\
224.000000000 0.832399004\\
256.000000000 1.020147478\\
};
\addlegendentry{$\np\!=\!4,\ \npj\!=\!1,\ \nv\!=\!4$}

\addplot[
    color={scalingcolor0},
    mark={*},
    mark options={solid},
    line width={0.75pt},
    mark size={0.8pt},
    dashed
]
table[row sep={\\}]{%
32.000000000 0.226033366\\
64.000000000 0.321751551\\
96.000000000 0.439364007\\
128.000000000 0.553211358\\
160.000000000 0.697483740\\
192.000000000 0.819679583\\
224.000000000 1.003750859\\
256.000000000 1.208076337\\
};
\addlegendentry{$\np\!=\!4,\ \npj\!=\!4,\ \nv\!=\!4$}
\addplot[
    color={scalingcolor2},
    mark={*},
    mark options={solid},
    line width={0.75pt},
    mark size={0.8pt},
    solid
]
table[row sep={\\}]{%
32.000000000 0.085972225\\
64.000000000 0.131436794\\
96.000000000 0.186944796\\
128.000000000 0.233347915\\
160.000000000 0.287006651\\
192.000000000 0.344014615\\
224.000000000 0.414194003\\
256.000000000 0.496664064\\
};
\addlegendentry{$\np\!=\!8,\ \npj\!=\!2,\ \nv\!=\!4$}

\addplot[
    color={scalingcolor2},
    mark={*},
    mark options={solid},
    line width={0.75pt},
    mark size={0.8pt},
    dashed
]
table[row sep={\\}]{%
32.000000000 0.228323989\\
64.000000000 0.269243163\\
96.000000000 0.336013682\\
128.000000000 0.380327484\\
160.000000000 0.436909787\\
192.000000000 0.499229949\\
224.000000000 0.595087881\\
256.000000000 0.662465696\\
};
\addlegendentry{$\np\!=\!8,\ \npj\!=\!8,\ \nv\!=\!4$}

\end{axis}
\end{tikzpicture}
        \caption{$M=12$, $J=4$ varying $N$ and different $(\np,\npj,\nv)$.}
        \label{fig:per_it_spring_mass_benchmark_p}
    \end{subfigure}
    \begin{subfigure}{\columnwidth}
        \centering
        \begin{tikzpicture}

\definecolor{scalingcolor0}{HTML}{1f77b4}

\begin{axis}[
    xlabel={Number of branches $J$},
    ylabel={Time per iteration [ms]},
    width={\textwidth},
    height={0.5\textwidth},
    ymin={0},
    ymax={0.6},
    grid={major},
    enlarge x limits={0.05},
    enlarge y limits={false},
    scaled ticks={false},
    tick align={outside},
    tick pos={left},
    major grid style={draw=black,draw opacity=0.15},
    legend pos={north west},
    legend cell align={left},
    legend style={
        font=\fontsize{5.5}{6.5}\selectfont,
        draw=black!20,
        fill=white,
        fill opacity=0.6,
        text opacity=1,
        inner xsep=2pt,
        inner ysep=1pt,
    },
    xtick={1,2,4,8}
]

\addplot[
    color={scalingcolor0},
    mark={*},
    mark options={solid},
    line width={0.75pt},
    mark size={0.8pt},
    solid
]
table[row sep={\\}]{%
1.000000000 0.083936466049\\
2.000000000 0.094981387755\\
4.000000000 0.131671263030\\
8.000000000 0.229011438667\\
};
\addlegendentry{$\np\!=\!8,\ \npj\!=\!\np/J,\ \nv\!=\!4$}

\addplot[
    color={scalingcolor0},
    mark={*},
    mark options={solid},
    line width={0.75pt},
    mark size={0.8pt},
    dashed
]
table[row sep={\\}]{%
1.000000000 0.078687061667\\
2.000000000 0.143930037778\\
4.000000000 0.267897828030\\
8.000000000 0.569295997333\\
};
\addlegendentry{$\np\!=\!8,\ \npj\!=\!8,\ \nv\!=\!4$}

\end{axis}

\end{tikzpicture}
        \caption{$M=12$, $N=64$ varying $J$.}
        \label{fig:per_it_spring_mass_benchmark_J}
    \end{subfigure}
    \par\vspace{0.5cm}
    \caption{Spring-Mass benchmark: Time per inner iteration of \FanCyQpalm{} with \texttt{parallel} (solid)  and \texttt{sequential} (dashed) branch execution.}    \label{fig:per_it_spring_mass_benchmark}
\end{figure}
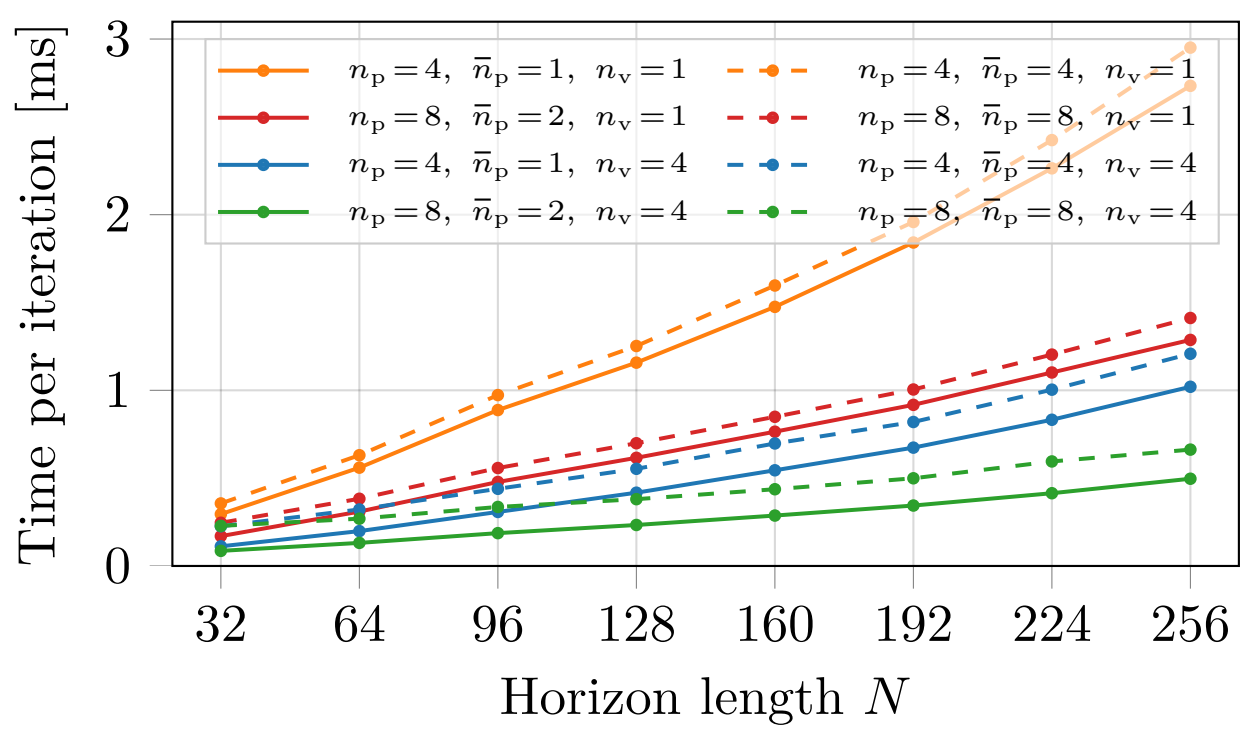
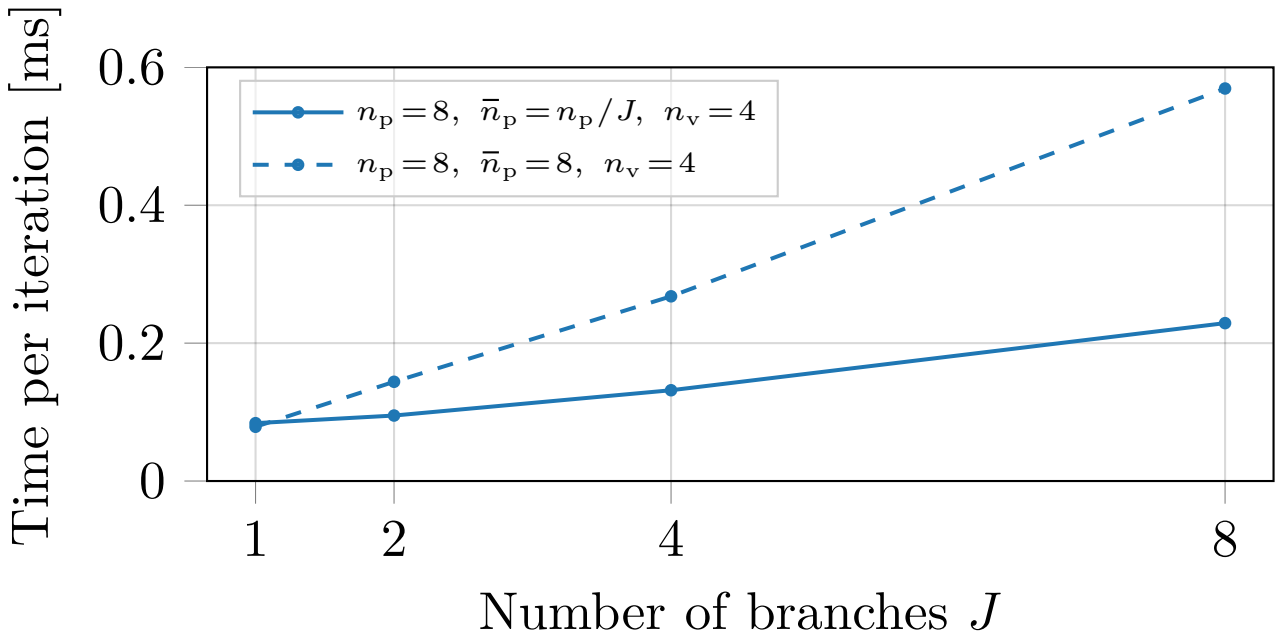

\begin{figure}[tbp]
    \centering
    \includegraphics[width=0.9\columnwidth]{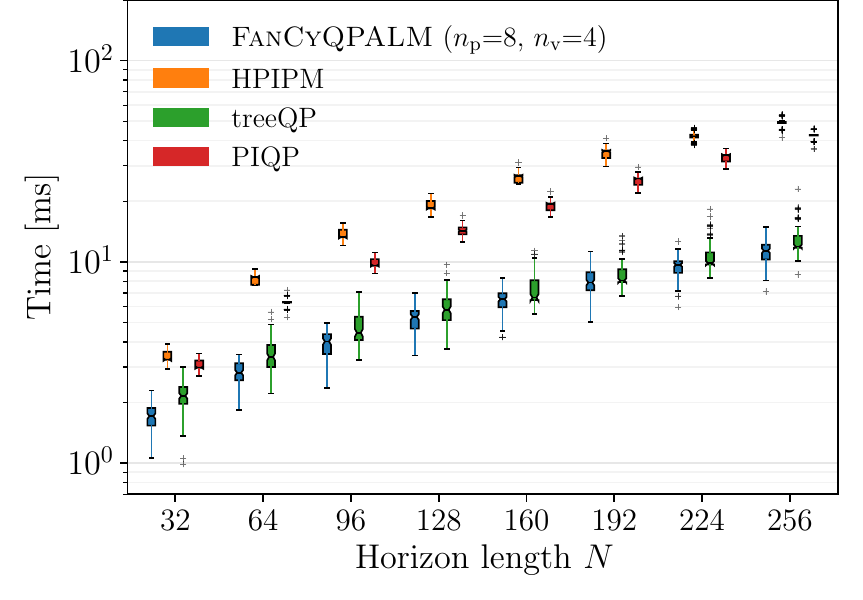}
    \caption{Spring-Mass benchmark: Total solver time with $M=12$ masses, $\numbranch=4$ branches and a varying horizon length $N$.}
    \label{fig:spring_mass_benchmark}
\end{figure}

\begin{figure}[tbp]
    \centering
    \includegraphics[width=0.9\columnwidth]{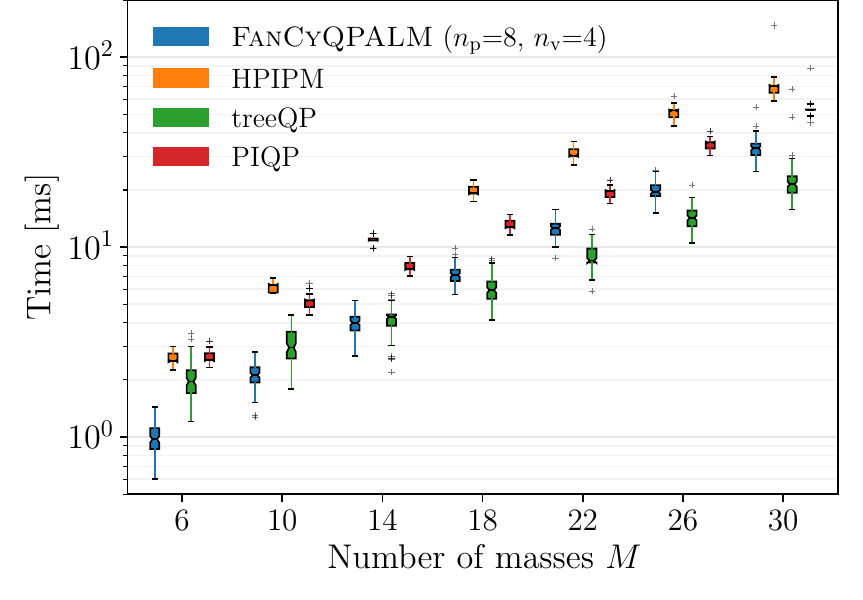}
    \caption{Spring-Mass benchmark: Total solver time with horizon length $N=64$, $\numbranch=4$ branches and a varying number of masses $M$.}
    \label{fig:spring_mass_benchmark_M}
\end{figure}

\Cref{fig:per_it_spring_mass_benchmark} shows the average time per inner iteration of the sequential and parallel variants of \FanCyQpalm{} for a spring-mass system with $M=12$ masses and varying horizon length $N$, number of branches $\numbranch$, total processor count $\np$, number of processors per branch $\npj$, and vector lengths $\nv$.
The parallel variant consistently outperforms the sequential variant (for $J>1$),
indicating that it is beneficial to parallelize over both branches and time.
For fixed $N,\np,$ and $\nv$, the performance gap increases as $J$ increases. For a fixed $J$ and varying $N$, the absolute performance gap between the sequential and parallel variants remains approximately constant.
\Cref{fig:spring_mass_benchmark,fig:spring_mass_benchmark_M}
show cold-start solve times for $\numbranch=4$, varying $N$
at $M=12$ and $M$ at $N=64$, respectively. Each boxplot
summarizes 100 instances with random initial states, using
the median solve time of ten cold-start solves per instance.

Across almost all tested instances, \FanCyQpalm{} achieves lower solve times than the fastest \textsc{Hpipm} and \PIQP{} configurations,
with the performance gap increasing with the horizon length. This is because \textsc{Hpipm} runs serially and therefore cannot make use of the additional compute power and cache available across multiple cores during factorization and solution.
In contrast, \PIQP{} requires two triangular solves per iteration.
For $M=12$, \FanCyQpalm{} achieves lower maximum times than \treeQP{},
despite the fact that (i) it takes more iterations, and (ii)
the problem involves only box constraints, which \treeQP{} exploits internally, 
whereas \FanCyQpalm{} treats them as general linear constraints with dense Jacobian matrices.

When varying $M$, \FanCyQpalm{} is superior for small $M$,
whereas \treeQP{} achieves the lowest median time for large $M$.
This crossover may be attributed to \treeQP{}'s repeated OpenMP
region entry and exit, compared with a single region for the
entire solve in \FanCyQpalm{}'s SPMD implementation.
Furthermore, \treeQP{}'s handling of box constraints involves 
a clipping operation, which is memory bound and scales linearly with M,
whereas the combined factorization and solve cost scales cubically.
Despite its efficiency, clipping can therefore be a
substantial share of runtime for small $M$, particularly during
the line search. As this overhead becomes less significant, \treeQP{}'s lower
iteration count can outweigh \FanCyQpalm{}'s lower per-iteration cost.

Each solution satisfies the prescribed tolerances of $10^{-8}$ for the
stationarity, equality, and inequality residuals and $10^{-6}$ for the
complementary-slackness residual. All residuals are evaluated externally
using the maximum ($\ell_\infty$) norm.%

\subsection{Quadruped Centroidal \ac{OCP} Benchmark} \label{sec:quadruped}
The quadruped centroidal \ac{OCP} formulation is based on the formulation in \cite[Section 8.2]{lezia2022}.
The contact sequence and foothold positions are predefined and the state vector $\ivar{x}{\bidx}{} \in \mathbb{R}^{12}$ is defined by stacking position, orientation, linear and angular velocities of the robot's center of mass.
The input vector $\ivar{u}{\bidx}{} \in \mathbb{R}^{12}$ is defined by stacking the contact forces at each of the four feet.
We consider the same parameters as in the reference implementation\footnote{https://gitlab.laas.fr/paleziart/quadruped-qp-problems}, but we regularize $Q\leftarrow Q+10^{-8}I$ as \treeQP{} requires positive definite cost matrices.
For the same reason, we take\footnote{The allowed range for $R$ is given in the reference implementation's README} $R= 5 \cdot 10^{-5}I$. The friction parameter $\mu_j$ is branch specific and is equidistantly varied between 0.3 and 0.7.
The temporal horizon $N\Delta t=0.64$ seconds stays fixed while $N$ is varied. We denote the legs by \ac{FL}, \ac{FR}, \ac{HL}, and \ac{HR}. The gait consists of four contact phases of equal duration (\ac{FL}-\ac{HR}, \ac{FR}-\ac{HL}, \ac{FL}-\ac{HR}, \ac{FR}-\ac{HL}), where each pair denotes the legs in contact with the ground. The longitudinal reference velocity remains constant at \(v_x=0.1\,\mathrm{m/s}\).
For \FanCyQpalm{}, we set $\texttt{initial\_inner\_tolerance}=10^{-4}$ instead of the default value 1, thereby solving the initial inner problems with higher accuracy, preventing excessively large penalty terms in the \ac{ALM}.
\begin{figure}[tbp]
\centering
    \includegraphics[width=0.9\columnwidth]{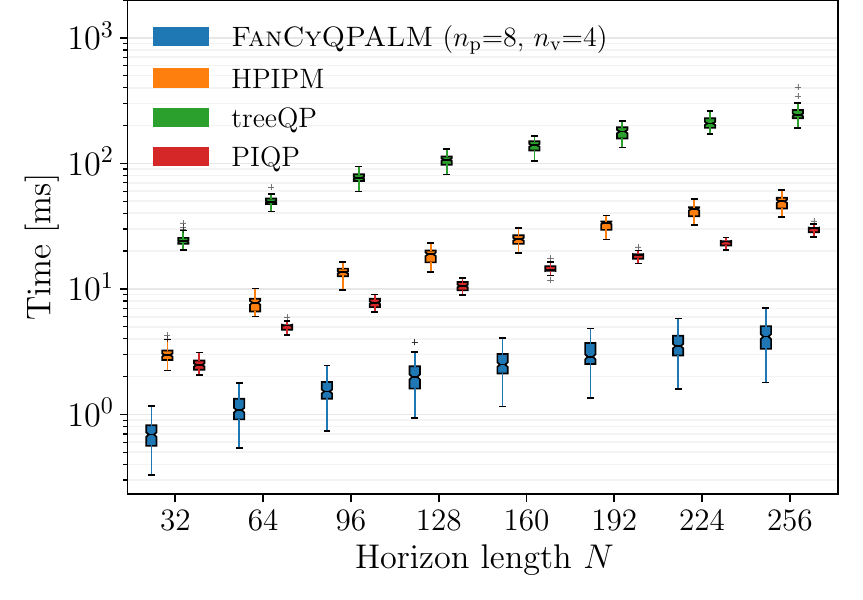}
    \caption{Quadruped benchmark: Total solver time for $\numbranch=4$ branches and a varying horizon length $N$.}
    \label{fig:quadruped_benchmark}
\end{figure}
\begin{figure}[tbp]
    \centering
    \includegraphics[width=0.9\columnwidth]{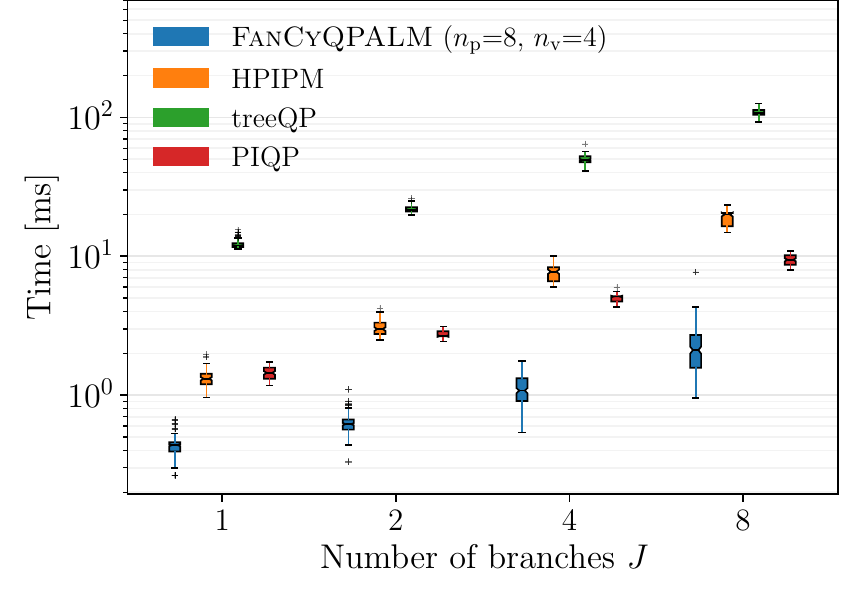}
    \caption{Quadruped benchmark: Total solver time for horizon $N\!=\!64$ and varying number of branches $\numbranch$.}
    \label{fig:quadruped_benchmark_J}
\end{figure}
In \Cref{fig:quadruped_benchmark}, we report the total cold-start solve time for the quadruped benchmark with $\numbranch=4$ and varying horizon length $N$.
In \Cref{fig:quadruped_benchmark_J}, we fix $N=64$ and vary the number of branches $\numbranch$.
 Each boxplot summarizes 100 instances with random initial positions $(x,y,z)$ uniformly distributed in $[-0.1,0.1]^2\times[0.135,0.335]$, using the median of ten solves per instance to reduce measurement noise.

\FanCyQpalm{} outperforms all other solvers for all tested horizon lengths and numbers of branches.
Most notably, in all tested configurations, treeQP's median solve time is more than 20 times that of \FanCyQpalm{}.
This is primarily due to the fact that the general stage-wise linear inequalities require \treeQP{} to solve its stage \acp{QP}
using qpOASES \cite{Ferreau2014}, rather than the more efficient \textit{clipping} approach mentioned above.

All solutions were checked externally. The stationarity, equality, and inequality residuals satisfy the tolerances of $10^{-8}$,
and the complementary-slackness residual satisfies the tolerance of $10^{-6}$, all measured in the $\ell_\infty$-norm.

\section{Conclusion}
We proposed \textsc{FanCyqlone}, a \ac{CPU}-based linear solver for multi-horizon optimal control problems that exploits parallelism across both the horizon and branches.
By extracting branch contributions directly from \textsc{Cyqlone} and combining them in a root system, the method efficiently solves multi-horizon KKT systems.
We integrated \textsc{FanCyqlone} into \QPALM{}, yielding \FanCyQpalm{} for linear \ac{OCP}s with general linear stage-wise constraints.
Benchmark comparisons with state-of-the-art solvers demonstrate the effectiveness of \FanCyQpalm{} on multi-horizon \ac{OCP}s.

\bibliographystyle{hieeetr}
\bibliography{references}

\begin{thebibliography}{10}

\bibitem{MPCbook}
J.~Rawlings, D.~Mayne, and M.~Diehl, {\em Model Predictive Control: Theory,
  Computation, and Design}.
\newblock Madison, WI, USA: Nob Hill Publishing, 2nd~ed., 2017.

\bibitem{lucia2013multi}
S.~Lucia, T.~Finkler, and S.~Engell, ``Multi-stage nonlinear model predictive
  control applied to a semi-batch polymerization reactor under uncertainty,''
  {\em Journal of Process Control}, vol.~23, no.~9, pp.~1306--1319, 2013.

\bibitem{gros20172stage}
E.~Klintberg and S.~Gros, ``A parallelizable interior point method for
  two-stage robust {MPC},'' {\em IEEE Transactions on Control Systems
  Technology}, vol.~25, no.~6, pp.~2087--2097, 2017.

\bibitem{multimodal2025}
S.~H. Nair, H.~Lee, E.~Joa, Y.~Wang, H.~E. Tseng, and F.~Borrelli, ``Predictive
  control for autonomous driving with uncertain, multimodal predictions,'' {\em
  IEEE Transactions on Control Systems Technology}, vol.~33, no.~4,
  pp.~1178--1192, 2025.

\bibitem{multimodel2025b}
K.~Ren, C.~Chen, H.~Sung, H.~Ahn, I.~M. Mitchell, and M.~Kamgarpour,
  ``Recursively feasible chance-constrained model predictive control under
  {Gaussian} mixture model uncertainty,'' {\em IEEE Transactions on Control
  Systems Technology}, vol.~33, no.~4, pp.~1193--1206, 2025.

\bibitem{9838190}
K.~S. Schweidel, S.~M. Koehler, V.~R. Desaraju, and M.~Barić,
  ``Driver-in-the-loop contingency {MPC} with invariant sets,'' in {\em 2022
  European Control Conference (ECC)}, pp.~808--813, 2022.

\bibitem{geurts2025contingency}
M.~Geurts, T.~Baltussen, A.~Katriniok, and M.~Heemels, ``A contingency model
  predictive control framework for safe learning,'' {\em IEEE Control Systems
  Letters}, vol.~9, pp.~1075--1080, 2025.

\bibitem{song2026accelerating}
F.~Song, L.~Zhang, L.~Wu, J.~Drgo{\v{n}}a, and C.~N. Jones, ``Accelerating
  branch {MPC} with two-level parallel direct solves on {GPUs},'' {\em arXiv
  preprint arXiv:2609.29977}, 2026.

\bibitem{Klintberg2016dual}
E.~Klintberg, J.~Dahl, J.~Fredriksson, and S.~Gros, ``An improved dual {Newton}
  strategy for scenario-tree {MPC},'' in {\em 2016 IEEE 55th Conference on
  Decision and Control (CDC)}, pp.~3675--3681, 2016.

\bibitem{alsterda2019contingency}
J.~P. Alsterda, M.~Brown, and J.~C. Gerdes, ``Contingency model predictive
  control for automated vehicles,'' in {\em 2019 American Control Conference
  (ACC)}, pp.~717--722, 2019.

\bibitem{alsterda2021contingency}
J.~P. Alsterda and J.~C. Gerdes, ``Contingency model predictive control for
  linear time-varying systems,'' {\em arXiv preprint arXiv:2102.12045}, 2021.

\bibitem{pas2025cyqlone}
P.~Pas and P.~Patrinos, ``Cyqlone: A parallel, high-performance linear solver
  for optimal control,'' {\em arXiv preprint arXiv:2512.09058v2}, 2025.

\bibitem{hermans2019qpalm}
B.~Hermans, A.~Themelis, and P.~Patrinos, ``{QPALM}: a {Newton-type} proximal
  augmented {Lagrangian} method for quadratic programs,'' in {\em 2019 IEEE
  58th Conference on Decision and Control (CDC)}, pp.~4325--4330, 2019.

\bibitem{hermans2022qpalm}
B.~Hermans, A.~Themelis, and P.~Patrinos, ``{QPALM}: a proximal augmented
  {Lagrangian} method for nonconvex quadratic programs,'' {\em Mathematical
  Programming Computation}, vol.~14, no.~3, pp.~497--541, 2022.

\bibitem{lowenstein2024qpalm}
K.~F. L{\o}wenstein, D.~Bernardini, and P.~Patrinos, ``{QPALM-OCP}: A
  {Newton-type} proximal augmented {Lagrangian} solver tailored for quadratic
  programs arising in model predictive control,'' {\em IEEE Control Systems
  Letters}, vol.~8, pp.~1349--1354, 2024.

\bibitem{frison2017high}
G.~Frison, D.~Kouzoupis, M.~Diehl, and J.~B. J{\o}rgensen, ``A high-performance
  {Riccati} based solver for tree-structured quadratic programs,'' in {\em 20th
  IFAC World Congress}, vol.~50, pp.~14399--14405, 2017.

\bibitem{frison2020hpipm}
G.~Frison and M.~Diehl, ``{HPIPM}: a high-performance quadratic programming
  framework for model predictive control,'' in {\em 21st IFAC World Congress},
  vol.~53, pp.~6563--6569, 2020.

\bibitem{Leidereiter2015dual}
C.~Leidereiter, A.~Potschka, and H.~G. Bock, ``Dual decomposition for {QPs} in
  scenario tree {NMPC},'' in {\em 2015 European Control Conference (ECC)},
  pp.~1608--1613, 2015.

\bibitem{kouzoupis2018dual}
D.~Kouzoupis, E.~Klintberg, M.~Diehl, and S.~Gros, ``A dual {Newton} strategy
  for scenario decomposition in robust multistage {MPC},'' {\em International
  Journal of Robust and Nonlinear Control}, vol.~28, no.~6, pp.~2340--2355,
  2018.

\bibitem{kouzoupis2019dual}
D.~Kouzoupis, E.~Klintberg, G.~Frison, S.~Gros, and M.~Diehl, ``A dual {Newton}
  strategy for tree-sparse quadratic programs and its implementation in the
  open-source software {treeQP},'' {\em International Journal of Robust and
  Nonlinear Control}, vol.~29, no.~8, pp.~2438--2457, 2019.

\bibitem{schwan2025exploiting}
R.~Schwan, D.~Kuhn, and C.~N. Jones, ``Exploiting multistage optimization
  structure in proximal solvers,'' in {\em 2025 IEEE 64th Conference on
  Decision and Control (CDC)}, pp.~4677--4683, 2025.

\bibitem{song2025parallel}
F.~Song, R.~Schwan, Y.~Chen, and C.~N. Jones, ``Parallel {KKT} solver in {PIQP}
  for multistage optimization,'' {\em arXiv preprint arXiv:2511.00946}, 2025.

\bibitem{sampathirao2024massively}
A.~K. Sampathirao, P.~Patrinos, A.~Bemporad, and P.~Sopasakis, ``Massively
  parallelizable proximal algorithms for large-scale stochastic optimal control
  problems,'' {\em Optimal Control Applications and Methods}, vol.~45, no.~1,
  pp.~45--63, 2024.

\bibitem{moran2025gpu}
R.~Moran and P.~Sopasakis, ``{GPU-accelerated SPOCK for scenario-based
  risk-averse optimal control problems},'' {\em arXiv preprint
  arXiv:2505.12078}, 2025.

\bibitem{bodard2023spock}
A.~Bodard, R.~Moran, M.~Schuurmans, P.~Patrinos, and P.~Sopasakis, ``{SPOCK}: A
  proximal method for multistage risk-averse optimal control problems,'' in
  {\em 22nd IFAC World Congress}, vol.~56, pp.~1944--1951, 2023.

\bibitem{zhang2025parallel}
L.~Zhang, C.~Lin, and S.~Grammatico, ``Parallel branch model predictive control
  on {GPUs},'' {\em arXiv preprint arXiv:2506.13624}, 2025.

\bibitem{sousa2026rake}
J.~Sousa-Pinto, ``Rake-compress {Riccati} recursions for parallel scenario-tree
  model predictive control,'' {\em arXiv preprint arXiv:2608.01332}, 2026.

\bibitem{BoydandWang2010}
Y.~Wang and S.~Boyd, ``Fast model predictive control using online
  optimization,'' {\em IEEE Transactions on Control Systems Technology},
  vol.~18, no.~2, pp.~267--278, 2010.

\bibitem{lezia2022}
P.-A. L\'{e}ziart, {\em Locomotion control of a lightweight quadruped robot}.
\newblock PhD thesis, Universit\'{e} Paul Sabatier - Toulouse III, 2022.

\bibitem{Ferreau2014}
H.~Ferreau, C.~Kirches, A.~Potschka, H.~Bock, and M.~Diehl, ``{qpOASES}: A
  parametric active-set algorithm for quadratic programming,'' {\em
  Mathematical Programming Computation}, vol.~6, no.~4, pp.~327--363, 2014.

\bibitem{crabtree1969identity}
D.~E. Crabtree and E.~V. Haynsworth, ``An identity for the {Schur} complement
  of a matrix,'' {\em Proceedings of the American Mathematical Society},
  vol.~22, no.~2, pp.~364--366, 1969.

\end{thebibliography}
\appendix

\begin{proof}[Proof of \cref{prop:factor}]
We omit the branch index $\bidx$ for notational simplicity.
The Riccati recursion and the \Cyqlone{} algorithm correspond to different elimination orders for the optimality conditions of the same augmented \ac{LQR} problem \eqref{eq:auglqr}.
After eliminating $\ivar{\tilde{x}}{}{0}$, the optimality conditions are given by the KKT system
$
    \mathscr K z = b$, or, as an augmented matrix
$
    \begin{pNiceArray}{c|c}[columns-width=1em]
    \mathscr K & b\\
    \end{pNiceArray}.
$

\Cyqlone{} computes a block Cholesky factorization of a permutation of the KKT
matrix $\mathscr K$ where the block row of the constraint $x^1_j = x$ and the
block column corresponding to its Lagrange multiplier are placed last.
That is, it considers the permuted system
\begin{equation*}
    \begin{pNiceArray}{c|c}[columns-width=1em]
        \mathscr K_\text{Cyq} & b_\text{Cyq}\\
    \end{pNiceArray}
    =
    \begin{pNiceArray}{cc|c}[columns-width=1em]
        \Pi_1 \mathscr K \Pi_1^\top &
        \Pi_1 \mathscr K \Pi_0^\top &
        \Pi_1 b \\
        \Pi_0 \mathscr K \Pi_1^\top &
        \Pi_0 \mathscr K \Pi_0^\top &
        \Pi_0 b \\
    \end{pNiceArray},
\end{equation*}
where $\Pi_1$ is the row selection matrix that represents the elimination order
of the \Cyqlone{} algorithm excluding the last $n_x$ rows, and $\Pi_0$ is the
row selection matrix that selects the last $n_x$ rows (therefore $\Pi_0 b=-x$).
By construction \cite[\S3.1]{pas2025cyqlone}, the bottom-right $n_x\times n_x$ block $L_0$ of the block Cholesky factor
computed by \Cyqlone{} is the Cholesky factor of the Schur complement
$
    \mathscr K_\text{Cyq}/(\Pi_1 \mathscr K \Pi_1^\top),
$
i.e., $L_0 L_0^\top = -\mathscr K_\text{Cyq}/(\Pi_1 \mathscr K \Pi_1^\top)$.
The bottom rows $\hat\lambda^0$ of the result of the forward
substitution step of \Cyqlone{} satisfy $
    -L_0 \hat\lambda^0 = \Pi_0 b - \Pi_0 \mathscr K \Pi_1^\top (\Pi_1 \mathscr K \Pi_1^\top)^{-1} \Pi_1 b,
$
or expressed as an augmented Schur complement system,
\begin{equation} \label{eq:schur-cyqlone}
    \begin{pNiceArray}{c|c}[columns-width=1em]
        \mathscr K_\text{Cyq} & b_\text{Cyq}\\
    \end{pNiceArray} / \bigl(\Pi_1 \mathscr K \Pi_1^\top\bigr)
    =
    \begin{pNiceArray}{c|c}[columns-width=1em]
        -L_0 L_0^\top & -L_0 \hat\lambda^0\\
    \end{pNiceArray}.
\end{equation}

Elimination of $\ivar{\tilde{x}}{}{0}$ followed by application of dynamic programming to the augmented \ac{LQR} problem \eqref{eq:auglqr}
yields the smaller problem
\begin{equation}
    \minimize_{\ivar{\tilde{u}}{}{0},\ivar{x}{}{1}} \; \frac{1}{2}\|\ivar{\tilde{u}}{}{0}\|^2 + V(\ivar{x}{}{1})\quad \sttlong\; \ivar{x}{}{1}=x,
\end{equation}
where $V$ is exactly the value function defined in \eqref{eq:lqr}. The resulting optimality conditions are given by
\begin{equation}
    \underbrace{\scalebox{0.9}{$\begin{pmatrix}
        \mathrm{I} & 0 & 0\\
        0 & P & -\mathrm{I} \\
        0 & -\mathrm{I} & 0 \\
    \end{pmatrix}$}}_{\mathscr K_\text{first}} \scalebox{0.9}{$\begin{pmatrix}
        \ivar{\tilde{u}}{}{0}\\
        \ivar{x}{}{1}\\
        \ivar{\mu}{}{}\\
    \end{pmatrix}$} = \underbrace{\scalebox{0.9}{$\begin{pmatrix}
        0\\
        -p\\
        -x\\
    \end{pmatrix}$}}_{b_\text{first}},
\end{equation}
Let
$\tilde\Pi_0$ be the row selection matrix that selects the rows corresponding to
the stationarity of $\tilde{u}^0$ and
$x^1$ and the constraint $x^1 = x$,
and let $\tilde\Pi_1$ select the remaining rows:
\begin{equation}
    \begin{pNiceArray}{c|c}[columns-width=1em]
        \mathscr K_\text{Ric} & b_\text{Ric}\\
    \end{pNiceArray} = \begin{pNiceArray}{cc|c}[columns-width=1em]
        \tilde\Pi_1 \mathscr K \tilde\Pi_1^\top &
        \tilde\Pi_1 \mathscr K \tilde\Pi_0^\top &
        \tilde\Pi_1 b \\
        \tilde\Pi_0 \mathscr K \tilde\Pi_1^\top &
        \tilde\Pi_0 \mathscr K \tilde\Pi_0^\top &
        \tilde\Pi_0 b \\
    \end{pNiceArray},
\end{equation}
Elimination using dynamic programming and Schur complements are equivalent \cite[\S9.2]{pas2025cyqlone},
so we have that
\begin{equation}
    \begin{pNiceArray}{c|c}[columns-width=1em]
        \mathscr K_\text{first} & b_\text{first}\\
    \end{pNiceArray}
    =
    \begin{pNiceArray}{c|c}[columns-width=1em]
        \mathscr K_\text{Ric} & b_\text{Ric}\\
    \end{pNiceArray} / (\tilde\Pi_1 \mathscr K \tilde\Pi_1^\top).
\end{equation}
Taking a final Schur complement in $\mathscr K_\text{first}$ 
to eliminate $\ivar{\tilde{u}}{}{0}$ and $\ivar{x}{}{1}$
yields the multiplier $\mu$ for the constraint $x^1 = x$ as the solution of the system
\begin{equation} \label{eq:schur-riccati}
    \begin{pNiceArray}{c|c}[columns-width=1em]
        \mathscr K_\text{first} & b_\text{first}\\
    \end{pNiceArray} / \begin{pmatrix}
        \mathrm I & 0\\
        0 & P
    \end{pmatrix} = \begin{pNiceArray}{c|c}[columns-width=1em]
        -P^{-1} & -P^{-1}p-x\\
    \end{pNiceArray}.
\end{equation}
The Schur complements \eqref{eq:schur-cyqlone} and \eqref{eq:schur-riccati} are equal because they eliminate the same rows and
columns of $\mathscr K$.
Successive and simultaneous elimination are equivalent
by the quotient rule for Schur complements \cite{crabtree1969identity}, and a permutation of the eliminated block does not change the Schur complement.
Therefore, $-P^{-1} = -L_0 L_0^\top$ and $-P^{-1}p-x = -L_0 \hat\lambda^0$. This implies that $L_0 = L^P$ in (\ref{F1})
(by uniqueness of the Cholesky factorization) and that $\hat\lambda^0 = \hat p$ in (\ref{S1}).
\end{proof}

\end{document}